\documentclass{cmmse2016}
\usepackage{amssymb}
\usepackage{amsmath,amsfonts,amssymb,amsthm,mathrsfs}
\usepackage[pdftex]{graphicx}
\usepackage{epstopdf}
\newtheorem{theorem}{Theorem}[]
\newtheorem{lemma}[theorem]{Lemma}

\theoremstyle{definition}

\begin{document}

%\pagestyle{empty}

%%%%%%%%%%%%%%%%%%
\def\R#1{\mathbb{R}^{#1}}
\def\B{\mathscr{B}}
\def\Lop{\mathcal{L}}
\def\Kop{\bar{\mathcal{K}}}
\def\opA#1{\mathcal{A}_{#1}}
\def\opB#1{\mathcal{B}_{#1}}
\def\half#1#2{\begin{matrix}\frac{#1}{#2}\end{matrix}}
\def\Char{\mathscr{C}}
\def\hate{e}
\def\Kopp{\widetilde{\mathcal{K}}}
\def\Rl{\widetilde{\mathscr{R}}_{\Lambda}}
\def\Rlo{\widetilde{\mathscr{R}}_0}
\def\thetamax{\Theta}

\newcommand{\Eqref}[1]{(\ref{#1})}

%  Data for the headings. Please fill both fields.

\markboth
  {Persistence analysis of  the age-structured population model} % abridged title
  {V. Kozlov, S. Radosavljevic, V. Tkachev, U. Wennergren}

%  Full title. Use '\\' to force a line break.

\title{Persistence  analysis of the age-structured population model on several patches}

%  Full author information.

%\author{Full name}{e-mail}{affiliation number (below)}

\author{Vladimir Kozlov}{vladimir.kozlov@liu.se}{1}

\author{Sonja Radosavljevic}{sonja.radosavljevic@liu.se}{1}

\author{Vladimir Tkachev}{vladimir.tkatjev@liu.se}{1}

\author{Uno Wennergren}{uno.wennergren@liu.se}{2}

%\affiliation{number}{department}{university}

\affiliation{1}{**Department of Mathematics}
  {Link\"oping University}

\affiliation{2}{**Department of Physics, Chemistry, and Biology}{Link\"oping University}

%  Abstract, key words and MSC codes.

\begin{abstract}
We consider a system of nonlinear partial differential equations that describes an age-structured population living in changing environment on $N$ patches. We prove existence and uniqueness of solution and analyze large time behavior of the system in time-independent case, for periodically changing and for irregularly varying environment. Under the assumption that every patch can be reached from every other patch, directly or through several intermediary patches, and that net reproductive operator has spectral radius larger than one, we prove that population is persistent on all patches. If the spectral radius is less or equal one, extinction on all patches is imminent.

\keywords Age-structured population model, temporal variability, spatial heterogeneity,  density-dependency, dispersion, large time behavior, periodic solutions

\msccodes AMS codes: 47N60, 47B65, 62P10
\end{abstract}

%
%  Main text of the article.
%

\section{Introduction}

Various natural processes and human activities are responsible for changes in habitat structure, quality and its fragmentation. A consequence of habitat reduction or destruction can be very severe and result in extinction of populations and loss of biodiversity. Understanding how habitat fragmentation and heterogeneity, climate change and seasonality interact with internal factors of population growth and influence population dynamics may provide a useful tool for conservation and managing of species (\cite{KW1}).

The idea that spatial structure is important and should be considered in studying population dynamics is not new. It is clear that there are no two identical habitats and that some are better than others. This can be linked to habitat heterogeneity and to the source-sink systems. A high quality habitat (source) yields positive growth rate of population, while low quality habitat (sink) has the opposite effect. For a species that inhabits several patches, possibility to move from one patch to another can be crucial for survival. Migration from source to sink has been studied and its stabilizing effect for population or even ecosystem is well-known. The importance of migration is in the fact that arrival of immigrants to the sink can save local population from extinction (\cite{PA}, \cite{PDias}).

No matter how important, spatial structure is not the only factor that influence population dynamics. It should be viewed as one of the pivotal factors and studied in combination with population structure, temporal variability and density-dependence (see for example \cite{BG}). All mentioned factors have relative importance for population dynamics which differs for terrestrial and marine species, vertebrate and invertebrate, plants and animals, large and small populations etc. The challenge is to make a model that includes all these factors and yet remains simple enough in order to be applicable to real world problems.

However, in spite of the evidence that various internal and external factors interact and contribute to population dynamics and complexity of ecological systems, a usual practice is to study only one or few of them, but not all of them simultaneously.  We mention some of them. Gurtin and McCamy introduced age-structured density dependent model in \cite{Grt}. Chipot considered age-structured  models with time and density-dependency in \cite{Chipot1} and \cite{Chipot2}.
Allen in \cite{Allen} studies different dispersion mechanisms for the logistic model without age-structure and without time-dependency in vital rates. Cui and Chen in \cite{CChen1} and \cite{CChen2} study the effect of diffusion on a single species or for predator-prey system, but they do not consider structured population. The system of delay differential equations developed by So, Wu and Zou in \cite{SWZ} deals with population divided into two age classes (immature and adult) and assumes that population inhabits two identical patches and that the vital rates are time-independent. In \cite{WXZ}, Weng, Xiao and Zou extended the model of So, Wu and Zou in \cite{SWZ} and considered dynamics of population on three patches. Terry  in \cite{Terry} investigates population of two stages and on two patches and discusses persistence of population for different birth functions and dispersion patterns. Takeuchi in \cite{Tak1} and \cite{Tak2} examined stability and effect of diffusion of generalized Lotka-Volterra systems. Hastings discussed stabilizing effect of dispersal in \cite{AH1}.

In contrast to the previous papers, we consider age-structure, time and density-depen\-dence and dispersion between patches and  investigate their impact on population dynamics. In formulation of the model we rely on several assumptions: population is age-structured and inhabits $N$ patches. In addition to this, environment changes with time, causing the change in the vital rates, competition level and dispersion coefficients. Competition for the resources occurs within each age-class which results in increased density-dependent mortality. The patches are not identical and they do not  have to be physically close, because many pest species exploit trading and transportation networks to move from one patch to another.

\section{Formulation of the model}

Under the made assumptions on mind, the system of balance equations for the model becomes:
\begin{align}\label{genpr}
\frac{\partial{n_k(a,t)}}{\partial t}+\frac{\partial{n_k(a,t)}}{\partial a} = -\mu_k(a,t)n_k(a,t)\left(1+\frac{n_k(a,t)}{L_k(a,t)}\right)  +\sum_{j=1}^ND_{kj}(a,t)n_j(a,t), \quad 1\le k\le N,
\end{align}
and the boundary and initial conditions are:
\begin{align}\label{genbc}
n_k(0,t)=\int_0^{\infty}m_k(a,t)n_k(a,t)\,da, \quad t> 0,
\end{align}
and
\begin{align}\label{genic}
n_k(a,0)=f_k(a), \quad a>0,
\end{align}
where $n_k(a,t)$ is the number of individuals of age $a$ at time $t$ on patch $k$, $\mu_k(a,t)$ is the death rate, $L_k(a,t)$ the regulating function, $$D=\|D_{kj}(a,t)\|_{1\le k,j\le N}
$$ the matrix of dispersion coefficients, $m_k(a,t)$ the birth rate and $f_k(a)$ the initial distribution of population.

The regulating function $L_k(a,t)$ represents limitations imposed on individuals by environment (or available resource per capita) and the logistic term $\frac{\mu_k(a,t)n_k^2(a,t)}{L_k(a,t)}$ describes increment in mortality due to population density. The regulation function resembles carrying capacity because for small values of
$L_k(a,t)$ mortality grows and for its large values, logistic term tends to zero.
Unlike carrying capacity it does not represent maximal supported population. We will use it as a starting point in study of relation between age-structure and density-dependence.

Underlying assumption is that population does not have to occupy all patches at initial time, but in order for population to survive, sufficiently young individuals must occupy at least one patch. Nonnegative dispersion coefficients $D_{kj}(a,t)$ for $k\ne j$ define proportion of individuals of age $a$ at time $t$ on patch $j$ that migrates to patch $k$. Every patch $k$ can be reached from every other patch $j$, possibly via several passing patches. Dispersion coefficients $D_{kk}(a,t)\le 0$ define proportion of population of age $a$ at time $t$ that leaves patch $k$.

Furthermore, we also assume  that the dispersion matrix is \textit{essentially positive}, i.e. for any $k\ne j$ there exist pairwise distinct indices $m_0, m_1,\ldots,m_s$ such that $m_0=k$, $m_s=j$ and $D_{m_i m_{i+1}}(a,t)>0$ for any $0\le i\le s-1$. The latter condition has the following natural explanation. Let us associate a directed graph $\Gamma(D,x)$ to the dispersion matrix $D$ as follows: the nodes of the graph are presented by the $N$ patches and the edges represent the transitions between patches according to the dispersion coefficients $D_{kj}(x)$ such that the $(k,j)$-entry of the incidence matrix of $\Gamma(D,x)$ is $0$ if $D_{kj}(x)=0$, and $1$ if $D_{kj}(x)>0$. Then $D$ satisfy the above condition  if and only if the graph $\Gamma(d,x)$ is connected (i.e. there is a path between every pair of vertices). An examples of an essentially positive  matrix is
$$
\left(
    \begin{array}{ccc}
      * & + & 0 \\
      0 & * & + \\
      + & 0 & * \\
    \end{array}
  \right)
$$
where $*$ means an entry of an unspecified sign.

%A life-history trade-off between reproduction and migration has been noted for many species, including migratory birds and some insects (see for example \cite{MZ}, \cite{CASW}, \cite{PAG}). This trade-off is caused by energy constraints because both reproduction and migration are energetically costly for organisms. Thus, the assumption that individuals do not reproduce during migration is biologically justified and mathematically it is stated as:
%\begin{align}\label{cond:D}
%\sum_{k=1}^ND_{kj}(a,t)\le |D_{jj}(a,t)|, \quad 1\le j\le N.
%\end{align}
%The relation between dispersion coefficients means that some migrants that are leaving patch $j$ can eventually die before reaching patch $k$, but they will not give birth during migration.

\section{The main results}
In this section we briefly describe our strategy and the principal results. The first main result is the existence and uniqueness of solution to the problem (\ref{genpr})--(\ref{genic}) in the class of bounded continuous functions. To this end, we consider the problem with time-independent vital rates, regulating function and dispersion coefficients. Namely, in order to determine the number of newborns
$$%\begin{align}\label{rho}
\rho_k(t):=n_k(0,t), \quad t\geq 0,\quad 1\leq k\leq N,
$$%\end{align}
we introduce two auxiliary initial value problems as follows. Let
$\Phi(x,y;\rho)$ and $\Psi(x,y;f)$  denote respectively the solutions of the following initial value problems:
\begin{equation*}
\begin{split}
  \frac{d h_k(x)}{dx}=-\mu_k(x,x+y) \left(h_k(x)+\frac{h_k^2(x)}{L_k(x,x+y)}\right)+\sum_{j=1}^ND_{kj}(x,x+y)h_j(x)=0, \,\,\,\, h(0)=\rho(y)
\\
  \frac{d h_k(x)}{dx}=-\mu_k(x+y,x) \left(h_k(x)+\frac{h_k^2(x)}{L_k(x+y,x)}\right) -\sum_{j=1}^ND_{kj}(x+y,x)h_j(x),
=0, \,\,\,\,
  h(0) = f(y).
\end{split}
\end{equation*}
It can be shown that each of the latter problems has a unique \textit{nonnegative} solution.

Then, the original problem can be reduced to the integral equation
\begin{align}\label{newb}
\rho(t)=\mathcal{K}\rho(t)+\mathcal{F}f(t),
\end{align}
where $\rho(t)=(\rho_1(t),\ldots,\rho_N(t))$ and the operators $\mathcal{K}$ and $\mathcal{F}$ are defined componentwise by
\begin{align}\label{k}
(\mathcal{K} \rho)_k(t)&=\int_0^tm_k(a,t)\Phi_k(a,t-a; \rho)\,da\\
\label{f}
(\mathcal{F} f)_k(t)&=\int_t^{\infty}m_k(a,t)\Psi_k(a,a-t; f)\,da,
\end{align}

In the time-independent case, an important role in description and analysis of solutions to (\ref{newb}) plays the so-called characteristic equation, i.e. the (unique) solution to
\begin{equation}\label{charact}
\rho=\Kop\rho.
\end{equation}
Here the operator $\Kop$ is given by
$$
(\Kop \rho)_k:=\int_0^{\infty}m_k(a)\Phi_k(a;\rho)\,da  , \qquad \rho\in \R{N}_+,
$$
i.e. when the newborns function $\rho$ is constant  for all $t\ge0$. Clearly, $\rho=0$ is always a solution of the characteristic equation, but our goal is to find out when a nontrivial solution $\rho\in \R{N}_+\setminus\{0\}$ does exist and is unique.

An important ingredient of our approach is the operator  ${\mathscr{R}}_0$ defined  by the right hand side in (\ref{genbc}), where $n_k(a,t)$ solves (\ref{genpr}) but without nonlinear term (formally assuming that $\frac{1}{L_k(a,t)}\equiv 0$) and with the boundary condition $n_k(0,t)=\rho_k=const$. We show that it can be alternatively defined as the blow-up of $\Kop$:
$$
{\mathscr{R}}_0\rho=\lim_{\varepsilon\rightarrow 0}\frac{1}{\varepsilon}\Kop\varepsilon \rho.
$$
Due to assumptions on dispersion coefficients $D_{kj}(a,t)$, this operator can be shown to be strongly positive and linear. Then an important corollary of the Krein-Rutman theorem is that its spectral radius $\sigma({\mathscr{R}}_0)$ is equal to the largest positive eigenvalue.

Thus defined operator $\mathscr{R}_0$ is called the \textit{net reproductive map} and $\sigma({\mathscr{R}}_0)$ is the net reproductive rate. Notice that in the one-dimensional case, $\sigma({\mathscr{R}}_0)$ coincides with the net reproductive rate $R_0$ which is  a well-established concept  and it is defined by
$$R_0=\int_0^{\infty}m(a)e^{-\int_0^a\mu(v)dv}\,da.$$
The net reproductive operator is just multiplication by this number.  In \cite{we2} we proved that $R_0$ corresponds to the characteristic equation
$$\int_0^{\infty}\frac{m(a)e^{-\int_0^a\mu(v)dv}}{1+n^*(1-e^{-\int_0^a\mu(v)dv})}\,da=1,$$
through which it is related to the average number of newborns $n^*$. Moreover, population declines for $R_0\le 1$ and grows for $R_0>1$.

We show that our definition is in consistence with the classical definition in one-dimensional case. Within this framework, we are able to characterize stationary solutions to the characteristic equation (\ref{charact}). Namely, the following dichotomy holds

\begin{theorem}\label{th1}
If $\sigma({\mathscr{R}}_0)\le 1$, then the characteristic equation has only trivial solution $\rho^*=0$. If $\sigma({\mathscr{R}}_0)> 1$, then the characteristic equation has exactly one nontrivial solution $\rho^*$ whose all components are positive.
\end{theorem}

One of ingredients of the proof is the following monotonicity result which has a considerable interest in itself. Here and in what follows we use the standard notation vector order relation: given two vectors $x,y\in \mathbb{R}^N$ one defines $x\le y$ if $x_i\le y_i$ for all $1\le i\le n$. Further, $x<y$ if $x\le y$ and $x\ne y$, and $x\ll y$ if $x_i< y_i$ for all $1\le i\le n$.

\begin{lemma}\label{lemma2}
Let $ w(x)=(w_1,...,w_N)$  be locally Lipschitz functions satisfying $w(0) \geq 0$ and
\begin{align}\label{i1}
w_k'(x)&\geq \sum_{j=1}^Nd_{kj}(x)w_j(x), \quad 1\leq k\leq N, \quad x\in [0,b),
\end{align}
where $d_{kj}(x)$ are continuous in $[0,b)$ and $d_{kj}(x)\geq 0$ for all $j\neq k$. Then $w(x)\geq 0$ on $[0,b)$. Furthermore, if $(d_{kj})_{1\le k,j\le N}$  is essentially positive and $w(0) > 0$ then $w(x) \gg 0$ on $(0,b)$.
\end{lemma}

Another important ingredient of the proof is a suitable generalization of the existence and the uniqueness results for monotone and concave operators established by  Krasnoselskii and Zabreiko in \cite{KrZa}.

%\subsection{Asymptotics of total population}
Our next result concerns solutions of the integral equation (\ref{newb}). We show that an a priori nonconstant solutions tends to the (constant) solution of the characteristic equation (\ref{charact}). In other words,

\begin{theorem}
Let $\rho=\rho(t)$ be a solution to (\ref{newb}). If $\sigma({\mathscr{R}}_0)\le 1$
then $\rho(t)\rightarrow 0$ as $t\rightarrow\infty$. If
$\sigma({\mathscr{R}}_0)> 1$, then $\rho(t)\rightarrow \rho^*$ as $t\rightarrow\infty$, where $\rho^*$ is defined by Theorem~\ref{th1}.
In particular,
$$\mbox{ if $\sigma({\mathscr{R}}_0)\le 1$, then $N(t)\rightarrow 0$ as $t\rightarrow\infty$,}$$
and
$$\mbox{ if  $\sigma({\mathscr{R}}_0)> 1$, then $N(t)\rightarrow \int_0^{\infty}\Phi(a;\rho^*)\,da$ as $t\rightarrow\infty$,}$$
where $N(t)=\int_0^{\infty}n(a,t)\,da$ is the total population and $\Phi(a;\rho^*)$ is a solution to the original problem with $\rho^*$ as the initial condition.
\end{theorem}

In this way, the net reproductive rate $\sigma({\mathscr{R}}_0)$ effectively determines large time behavior of population on $N$ patches in constant environment. Here, as in the one-dimensional case, $\sigma({\mathscr{R}}_0)\le 1$ implies extinction of population on all patches, while $\sigma({\mathscr{R}}_0)> 1$ grants persistence of population.

\section{Periodically varying environment}

Most natural habitats are positively autocorrelated, see for example  \cite{Steele}. Thus, the assumption that the vital rates, regulating function and dispersal coefficient are changing periodically with respect to time is reasonable. In studying large time behavior of solution to equation (\ref{newb}) in periodically changing environment, a pivotal role belongs to the characteristic equation $\rho(t)=\Kopp\rho(t)$, where operator $\Kopp$ is given by the right hand side of (\ref{genbc}) and $n_k(a,t)$ solves (\ref{genpr}) with periodic boundary condition $n_k(0,t)=\rho_k(t)$, $t>0$.

Similarly to the previous situation, we introduce the net reproductive operator $\Rlo$ by
the right hand side of (\ref{genbc}) assuming that $n_k(a,t)$ solves (\ref{genpr}) without nonlinear term and with periodic boundary condition $n_k(0,t)=\rho_k(t)$. In this case, we have that
$$
\Rlo \rho =\lim_{\varepsilon\rightarrow 0}\frac{1}{\varepsilon}\Kopp\varepsilon \rho.
$$
Operator $\Rlo$ is strictly positive, linear and defined on space of periodic continuous functions. Its spectral radius $\sigma(\Rlo)$ is equal to the largest eigenvalue and it is called the net reproductive rate.

One of the main results about periodic case is that if $\sigma(\Rlo)\le 1$, then the characteristic equation $\rho(t)=\Kopp\rho(t)$ has only trivial solution $\rho\equiv 0$. If $\sigma(\Rlo)>1$, then the characteristic equation has exactly one nontrivial solution $\rho^*(t)$, where all components are positive periodic functions. Furthermore, we show that the number of newborns in periodic case converges to a unique periodic solution of the characteristic equation.
For the total population, the result can be formulated in the following way:
$$\mbox{ if $\sigma(\Rlo)\le 1$, then $N(t)\rightarrow 0$ as $t\rightarrow\infty$,}$$
and
$$\mbox{ if  $\sigma(\Rlo)> 1$, then $N(t)\rightarrow \int_0^{\infty}\Phi(a,t-a;\rho^*)\,da$ as $t\rightarrow\infty$,}$$
where $\Phi(a,t;\rho^*)$ is a solution to the original initial value problem and $\rho^*(t)$ is a periodic solution to the characteristic equation. In this situation, as in the time-independent case, the net reproductive rate $\sigma(\Rlo)$ determines extinction or persistence of population on all patches.

Upper and lower bounds for population growth can be found even if environment is changing irregularly. Namely, if the parameters of the original model are estimated from above and below by periodic functions for large time, then two periodic problems can be formulated. One of these periodic problems is an upper bound for the original problem, and the other is a lower bound. The number of newborns in the original problem, $\rho(t)$, is then estimated by the number of newborns in periodic problems. In other words, the following holds:
$$
\rho^-(t)-\varepsilon\le\rho(t)\le\rho^+(t)+\varepsilon,
$$
where $\rho^{\pm}(t)$ are solutions of the characteristic equations for periodic problems and $\varepsilon>0$ is an arbitrary small number.

In order to explain the influence of dispersion on persistence of population, we compare the system of $N$ patches with dispersion with the system of $N$ isolated patches. To obtain the net reproductive operator $\mathscr{R}_{0}$ and the net reproductive rate $\sigma(\mathscr{R}_{0})$, we use system (\ref{genpr}) without nonlinear term. In the worst case scenario for the system with dispersion, all migrants die before reaching next patch i.e. $D_{kj}(a)=0$ for $a>0$, $1\le k,j\le N$, $k\ne j$. In this case we estimate $n_k(a,t)$ in the following way
$$\frac{\partial{n_k(a,t)}}{\partial t}+\frac{\partial{n_k(a,t)}}{\partial a} \ge -(\mu_k(a)+|D_{kk}(a)|)n_k(a,t)$$
and obtain
$$n_k(a,t)\ge \rho_ke^{-\int_
0^a(\mu_k(v)+|D_{kk}(v)|)dv}\,da.$$
This implies that
$$\sigma(\mathscr{R}_{0})\rho_k=\int_0^{\infty}m_k(a)n_k(a,t)\,da \ge \rho_k \int_0^{\infty}m_k(a)e^{-\int_0^a(\mu_k(v)+|D_{kk}(v)|)dv}\,da$$
and
\begin{align}\label{est:sigma}
\sigma(\mathscr{R}_{0})\ge \max_{1\le k\le N}\int_0^{\infty}m_k(a)e^{-\int_0^a(\mu_k(v)+|D_{kk}(v)|)dv}\,da
\end{align}
On the other hand, the system of isolated patches corresponds to the problem (\ref{genpr})--(\ref{genic}) with $D_{kj}(a,t)=0$ for $a,t >0$ and $1\le k,j \le N$. The net reproductive rate on each patch is
$$\sigma_k=\int_0^{\infty}m_k(a)e^{-\int_0^a\mu_k(v)dv}\,da.$$
If $\sigma_k\le 1$ for some $k$, extinction of population on patch $k$ is imminent. However, according to (\ref{est:sigma}), it follows that $\sigma(\mathscr{R}_{0})>1$ if there exists at least one patch such that $\sigma_k>1$ and $|D_{kk}(a)|$ is small enough. This means that persistence of population on all patches is possible under assumption that population would persist on at least one patch and that the number of individuals emigrating from this patch is sufficiently small.

\section{Irregularly varying environment}
Previously, we have established conditions that ensures persistence of population in time-independent case and for periodically changing environment.
Here we establish  asymptotic behavior of the solution to the main model in the general time-dependent case in two ways: by using periodic functions to formulate upper and lower bounds and
by finding upper and lower boundaries for the sum of newborns on all patches and for the total population by relaxing conditions for dispersion.  In the first case, we will find conditions for extinction or persistence of population under assumption that the vital rates, regulating function and dispersion coefficients are bounded by periodic functions. In the second case we obtain a crude estimate of newborns density and total population which will not depend on migration pattern.

Suppose that the vital rates, regulating function and dispersion coefficients have periodic bounds, or, in other words, that the following estimates hold for the structure coefficients and large~$t$:
\begin{equation}\label{pb}
\begin{aligned}
m^-(a,t) &\le m(a,t) \le m^+(a,t), \\
\mu^+(a,t) &\le \mu(a,t) \le \mu^-(a,t), \text{ etc}
\end{aligned}
\end{equation}
where  $m^{\pm}_k$, $\mu^{\pm}_k$, $L^{\pm}_k$ and $D^{\pm}_{kj}$ are some $T$-periodic functions. Instead of the original problem (\ref{genpr}), we study two associated periodic problems with parameters defined by (\ref{pb}). It is natural to expect that $n(0,t)$ will be bounded from above and below by $\rho^{\pm}(t)$ for sufficiently large $t$, where $n(0,t)$ is a number of newborns in the original problem and $\rho^{\pm}(a,t)$ are solutions to the characteristic equations related to periodic problems
\begin{align*}
\Kopp^{\pm}\rho(t)=\rho(t), \quad t\in\R{},
\end{align*}
where operators $\Kopp^{\pm}$ are defined componentwise by
\begin{align}\label{charpKpm}
(\Kopp^{\pm}\rho)_k(t):=\int_0^{\infty}m_k^{\pm}(a,t)\Phi_k^{\pm}(a,t-a;\rho)\,da, \quad t\in\mathbb{R}, \quad 1\leq k\leq N.
\end{align}
The functions $\Phi^{\pm}(x,y;\rho)$ are unique solutions to the initial value problems
\begin{equation}\label{PDEpm}
\renewcommand\arraystretch{1.5}
\left\{
\begin{array}{rcl}
\frac{d}{dx}h_k^{\pm}(x)&=&-\mu_k^{\pm}(x,x+y) \left(1+\frac{h_k^{\pm}(x)}{L_k^{\pm}(x,x+y)}\right)h_k^{\pm}(x) +\sum_{j=1}^ND_{kj}^{\pm}(x,x+y)h_j^{\pm}(x),
 \\
  h^{\pm}(0) &=& \rho(y),
\end{array}
\right.
\end{equation}
where the coefficients are T-periodic and satisfy condition (\ref{pb}) and the initial conditions are given by a vector-function $\rho\in C(\R{}_+,\R{N}_+)$ such that $\rho(t+T)=\rho(t), \quad t\in \R{}.$ Then
\begin{align*}
\rho^{\pm}(t)=\mathcal{K}^{\pm}\rho^{\pm}(t)+\mathcal{F}^{\pm}f(t),
\end{align*}
where operators $\mathcal{K}^{\pm}$ and $\mathcal{F}^{\pm}$ are defined component-wise by
\begin{align*}
(\mathcal{K}^{\pm}\rho)_k(t)=\int_0^tm^{\pm}_k(a,t)\Phi^{\pm}_k(a,t-a;\rho)\,da,\\
(\mathcal{F}^{\pm}f)_k(t)=\int_t^{\infty}m^{\pm}_k(a,t)\Psi^{\pm}_k(a,a-t;f)\,da.
\end{align*}

The corresponding net reproductive operators and net reproductive rates are denoted by $\Rlo^{\pm}$ and $\sigma(\Rlo^{\pm})$. Then the next result states  that the number of newborns in irregularly changing environment can be bounded from above and below by the number of newborns in the associated periodically changing environments.

\begin{theorem}
Let $\rho(t)$ be a solution to equation (\ref{newb}).
If $\sigma(\Rlo^+)\le 1$, then $\rho(t)\rightarrow 0$ as $t\rightarrow\infty$. If $\sigma(\Rlo^-)>1$, then
\begin{align}\label{est:pernb}
\rho^-(t)-\varepsilon\le \rho(t) \le \rho^+(t)+\varepsilon \quad\mbox{for large $t$},
\end{align}
where $\rho^{\pm}(t)$ are solutions to (\ref{charpKpm}) and $\varepsilon$ is an arbitrary positive number.
\end{theorem}

%
%  Bibliography. Follow the usual conventions.
%
\def\cprime{$'$} \def\cprime{$'$}
\providecommand{\bysame}{\leavevmode\hbox to3em{\hrulefill}\thinspace}
\providecommand{\MR}{\relax\ifhmode\unskip\space\fi MR }
\providecommand{\MRhref}[2]{%
  \href{http://www.ams.org/mathscinet-getitem?mr=#1}{#2}
}
\providecommand{\href}[2]{#2}


\begin{thebibliography}{10}

\bibitem{Allen}
L.~Allen, \emph{Persistence and extinction in single-species reaction-diffusion
  models}, Bulletin of Mathematical Biology \textbf{45} (1983), 209--227.

\bibitem{PA}
P.~Amarasekare, \emph{The role of density-dependent dispersal in source-sink
  dynamics}, Journal of Theoretical Biology \textbf{226} (2004), 159–168.

\bibitem{BG}
O.~N. Bjørnstad and B.~T. Grenfell, \emph{Noisy clockwork: Time series
  analysis of population fluctuations in animals}, Science Translational
  Medicine \textbf{293 (5530)} (2001), 638--643.

\bibitem{CChen1}
J.~Chi and L.~Chen, \emph{The effect of diffusion on the time varying logistic
  population growth}, Computers Math. Applic. \textbf{36} (1998), 1--9.

\bibitem{CChen2}
\bysame, \emph{Permanence and extinction in logistic and lotka-volterra systems
  with diffusion}, Journal of Mathematical Analysis and Applications
  \textbf{258} (2001), 512--535.

\bibitem{Chipot1}
M.~Chipot, \emph{On the equations of age-dependent population dynamics}, Arch.
  Rational Mech. Anal. \textbf{82} (1983), no.~1, 13--25.

\bibitem{Chipot2}
\bysame, \emph{A remark on the equation of age-dependent population dynamics},
  Quarterly of Applied Mathematics \textbf{42} (1984), no.~2, 221--224.

\bibitem{PDias}
P.~C. Dias, \emph{Sources and sinks in population biology}, Trends in Ecology
  and Evolution \textbf{11} (1996), 326--330.

\bibitem{CASW}
C.~A. Schmidt-Wellenburg et. al., \emph{Trade-off between migration and
  reproduction: does a high workload affect body condition and reproductive
  state?}, Behavioral Ecology \textbf{doi:10.1093/beheco/arn066} (2008).

\bibitem{PAG}
P.~A. Guerra, \emph{Evaluating the life-history trade-off between dispersal
  capability and reproduction in wing dimorphic insects: a meta-analysis},
  Biological Reviews \textbf{86} (2011), 813–835.

\bibitem{Grt}
M.~E. Gurtin and R.~C. MacCamy, \emph{Nonlinear age-dependent population
  dynamics}, Arch. Rat. Mech. Anal. \textbf{54} (1974), 281--300.

\bibitem{AH1}
A.~Hastings, \emph{Complex interactions between dispersal and dynamics: Lessons
  from coupled logistic equations}, Ecology \textbf{44} (1993), 1362--1372.

\bibitem{KW1}
P.~Kareiva and U.~Wennergren, \emph{Connecting landscape patterns to ecosystem
  and population processes}, Nature \textbf{373} (1995), 299--302.

\bibitem{KrZa}
M.~A. Krasnosel{\cprime}ski{\u\i} and P.~P. Zabre{\u\i}ko, \emph{Geometrical
  methods of nonlinear analysis}, Grundlehren der Mathematischen Wissenschaften
  [Fundamental Principles of Mathematical Sciences], vol. 263, Springer-Verlag,
  Berlin, 1984, Translated from the Russian by Christian C. Fenske. \MR{736839}

\bibitem{MZ}
S.~Mole and A.~Zera, \emph{Differential allocation of resources underlies the
  dispersal-reproduction trade-off in the wing-dimorphic cricket, gryllus
  rubens}, Oecologia \textbf{93} (1993), 121--127.

\bibitem{SWZ}
J.~W.-H. So, J.~Wu, and X.~Zou, \emph{Structured population on two patches:
  modeling dispersal and delay}, J. Math. Biol. \textbf{43} (2001), 37--51.

\bibitem{Steele}
J.~H. Steele, \emph{A comparison of terrestrial and marine ecological systems},
  Nature \textbf{313} (1985), 355--358.

\bibitem{Tak1}
Y.~Takeuchi, \emph{Diffusion effect on stability of lotka-volterra models},
  Bul. Math. Bio. \textbf{48} (1986), 585--601.

\bibitem{Tak2}
\bysame, \emph{Global stability in generalized lotka-volterra diffusion
  systems}, Journal of Math. Analysis and Appl. \textbf{116} (1986), 209--221.

\bibitem{Terry}
A.~Terry, \emph{Dynamics of structured population on two patches}, Journal of
  Mathematical Analysis and Applications \textbf{378} (2011), 1--15.

\bibitem{WXZ}
P.~Weng, C.~Xiao, and X.~Zou, \emph{Rich dynamics in a non-local population
  model over three patches}, Nonlinear Dyn. \textbf{59} (2010), 161--172.

\bibitem{we2}
V.~Kozlov, S~Radosavljevic, and U.~Wennergren, \emph{Large time behavior of
  logistic age-structured population model in changing environment}, submitted.



\end{thebibliography}
\end{document}